\documentclass[a4paper,draft]{amsproc}
\usepackage{amssymb}
\usepackage[hyphens]{url} 

\theoremstyle{plain}
 \newtheorem{thm}{Theorem}[section]

\theoremstyle{definition}

\theoremstyle{remark}
 
 \numberwithin{equation}{section}

\renewcommand{\leq}{\leqslant}
\renewcommand{\geq}{\geqslant}

\title[A matrix application of quasi-monotone sequences to Fourier series]{A matrix application of quasi-monotone sequences to Fourier series}



\author[YILDIZ]{\bfseries \c{S}ebnem Y{\i}ld{\i}z} 

\address{Department of Mathematics \\
Ahi Evran University\\ K{\i}r\c{s}ehir, Turkey}
\email{sebnemyildiz@ahievran.edu.tr; sebnem.yildiz82@gmail.com}

\begin{document}

{\begin{flushleft}\baselineskip9pt\scriptsize
\end{flushleft}}
\vspace{18mm} \setcounter{page}{1} \thispagestyle{empty}

\begin{abstract}
In this paper, we have generalized a main theorem dealing with weighted mean summability method for absolute matrix summability method which plays a vital role in  summability theory and applications to the other sciences by using quasi-monotone sequences. The main result in this paper extends the results in \cite{Bor4}.

\end{abstract}

\maketitle

\section{Introduction}
A sequence $(d_{n})$ is said to be $\delta$-quasi-monotone, if $d_{n}\rightarrow0$, $d_{n}>0$ ultimately, and $\Delta d_{n}\geq-\delta_{n}$, where  $\Delta d_n = d_n -d_{n+1}$ and $\delta$=$(\delta_{n})$ is a sequence of positive numbers (see \cite{Boas}). For any sequence $(\lambda_{n})$ we write that $\Delta\lambda_n=\lambda_n-\lambda_{n+1}$.
The sequence $(\lambda_{n})$ is said to be of bounded variation,
denoted by $(\lambda_{n})\in  BV$, if
$\sum_{n=1}^{\infty}\left\vert \Delta\lambda_{n}\right\vert <\infty$. Let $\sum a_{n}$ be a given infinite series with partial sums $(s_{n})$. We denote by $u_n^\alpha$ and $t_n^\alpha$ the $n$th Ces\`aro means of order $\alpha$, with $\alpha>-1$, of the
sequences $(s_n)$ and $(na_n)$, respectively, that is (see \cite{Ce}),

\begin{equation}
u_{n}^{\alpha}=\frac{1}{A_{n}^{\alpha}}\sum_{v=0}^{n}A_{n-v}^{\alpha-1}s_{v}~~~~
{\texttt{and}}~~~~
t_{n}^{\alpha}=\frac{1}{A_{n}^{\alpha}}\sum_{v=1}^{n}A_{n-v}^{\alpha-1}va_{v}, ~~ ({t_{n}}^{1}=t_{n})
\end{equation}
\noindent where
\begin{equation} A_{n}^{\alpha}=\frac{(\alpha+1)(\alpha+2)....(\alpha+n) }{n!}=
O(n^\alpha), \quad {A_{-n}^\alpha=0} \quad {\texttt{for}} \quad  n>0.\
\end{equation}
A series $\sum{a_n}$ is said to be summable
$\mid{C},\alpha\mid_k$, $k\geq1$, if (see \cite{Fl}, \cite{Ko})
\begin{equation}
\sum_{n=1}^{\infty} n^{k-1}\mid
u_n^{\alpha}-u_{n-1}^{\alpha}\mid^{k} = \sum_{n=1}^{\infty} \frac{1}{{n}}\mid t_n^{\alpha}\mid^k <
\infty.
\end{equation}
If we set $\alpha$=1, then we have ${\mid{C,1}\mid}_k$ summability.
Let $(p_{n})$ be a sequence of positive number such that
\begin{equation}
P_{n}=\sum_{v=0}^{\infty}p_{v}\rightarrow \infty \quad {as}\quad
{n}\rightarrow \infty ,\quad (P_{-i}=p_{-i}=0,~~i\geq1).
\end{equation}\\

The sequence-to-sequence transformation
\begin{equation}
w_{n}=\frac{1}{P_{n}} \sum_{v=0}^{n}p_{v}s_{v}
\end{equation}
defines the sequence $(w_{n})$ of the Riesz mean or simply the $\left(\bar{N},p_{n}\right)$
mean of the sequence $(s_{n})$, generated by the sequence of coefficients $(p_{n})$ (see \cite{Ha}).
The series $\sum{a_n}$ is said to be summable
$|\bar{N},p_{n}|_{k}$, $k\geq1$, if (see \cite{Bor1})
\begin{equation}
\sum_{n=1}^{\infty}\left(\frac{P_{n}}{p_{n}}\right)^{k-1}\left|w_{n}-w_{n-1}\right|^{k} < \infty.
\end{equation}
In the special case when $p_{n}=1$ for all values of $n$ (respect. $k=1$), then $|\bar{N},p_{n}|_{k}$ summability is the same as $|C,1|_{k}$ (respect. $|\bar{N},p_{n}|$) summability. We write
$X_{n}=\sum_{v=1}^{n}\frac{p_{v}}{P_{v}}$, then $(X_{n})$ is a positive increasing sequence tending to infinity with $n$.\\
\subsection{ An application to trigonometric Fourier series}
Let $f$ be a periodic function with period $2\pi$ and integrable $(L)$ over $(-\pi,\pi)$. Without any loss of generality  the constant term in the Fourier series of $f$ can be taken to be zero, so that
\begin{equation}
f(t)\sim \sum_{n=1}^{\infty}(a_{n}cosnt+b_{n}sinnt)=\sum_{n=1}^{\infty}C_{n}(t).
\end{equation}
where
\begin{equation}
a_{0}=\frac{1}{\pi}\int_{-\pi}^{\pi}f(t)dt,\quad\nonumber
a_{n}=\frac{1}{\pi}\int_{-\pi}^{\pi}f(t)cos(nt)dt, \nonumber\quad
b_{n}=\frac{1}{\pi}\int_{-\pi}^{\pi}f(t)sin(nt)dt.\nonumber
\end{equation}
We write
\begin{equation}
\phi(t)=\frac{1}{2}\left\lbrace f(x+t)+f(x-t)\right\rbrace,\\
\phi_{\alpha}(t)= \frac{\alpha}{t^{\alpha}} \int_{0}^{t} (t-u)^{\alpha-1}\phi(u)\, du, \quad (\alpha> 0).
\end{equation}
It is well known that if $\phi(t)\in \mathcal{BV}(0,\pi)$, then $z_{n}(x)=O(1)$, where  $z_{n}(x)$ is the $(C,1)$ mean of the sequence $(nC_{n}(x))$ (see \cite{Ch}).\\
\begin{equation}
z_{n}(x)=\frac{1}{n+1}\sum_{v=1}^{n}vC_{v}(x).
\end{equation}
Let $A=(a_{nv})$ be a normal matrix,
i.e., a lower triangular matrix of nonzero diagonal entries. Then $A$ defines the sequence-to-sequence transformation, mapping the
sequence $s=(s_{n})$ to $As=\left(A_{n}(s)\right)$, where
\begin{equation}\label{eq:6}
A_{n}(s)=\sum_{v=0}^{n}a_{nv}s_{v}, \quad n=0,1,...
\end{equation}
The series $\sum a_{n}$ is said to be summable $\left|A,p_{n} \right|_{k}$, $k\geq 1$, if (see \cite{Su})
\begin{equation}\label{eq:7}
\sum_{n=1}^{\infty}\left(\frac{P_{n}}{p_{n}} \right) ^{k-1}\left|\bar{\Delta}A_{n}(s)\right|^{k}< \infty,
\end{equation}
where
\begin{equation}
\bar{\Delta}A_{n}(s)=A_{n}(s)-A_{n-1}(s).
\end{equation}
If we take $p_{n}=1$, for all $n$, $\left|A,p_{n} \right|_{k}$ summability is the same as  $\left| A \right|_{k}$ summability (see \cite{NT}). And also if we take  $a_{nv}=\frac{p_{v}}{P_{n}}$, then we have $\left| \bar{N},p_n\right|_k$ summability. \\

\begin{thm}   \cite{Bor4}
	Let $\lambda_{n}\rightarrow0$ as $n\rightarrow\infty$ and let $(p_{n})$ be a sequence of positive numbers such that
	\begin{equation}
	P_{n}=O(np_{n}) \ ~~ as \ n\rightarrow\infty.
	\end{equation}
	Suppose that there exists a sequence of numbers $(A_n)$ which is $\delta$-quasi-monotone with $\sum n X_n \delta_n<\infty$, $\sum A_n X_n$ is convergent, and $|\Delta\lambda_n|\leq|A_n|$ for all $n$. If
	\begin{equation}
	\sum_{n=1}^{m}\frac{p_{n}}{P_{n}}\frac{|t_{n}|^{k}}{{X_{n}}^{k-1}}=O(X_{m})\quad
	as \quad m\rightarrow\infty,
	\end{equation}
	satisfies, then the series $\sum a_{n}\lambda_{n}$ is summable $|\bar{N},p_{n}|_{k}$, $k\geq1$.\\
\end{thm}

Bor has obtained the following result dealing with Fourier series.\\
\begin{thm} \label{TheoremBor2} \cite{Bor4}  If $\phi_1(t)\in \mathcal{BV}(0,\pi)$, and the
	sequences $(A_n)$, $(\lambda_n)$, and $(X_n)$ satisfy the conditions
	of Theorem 1.1, then the series  $\sum C_{n}(x)\lambda_n $ is
	summable $|\bar{N},p_{n}|_{k}$, $k\geq1$.
\end{thm}
\section {MAIN RESULTS}
The aim of this paper is to generalize Theorem 1.2 for $|A,p_{n}|_{k}$ summability factors of Fourier series.
\\Given a normal matrix $A=(a_{nv})$, we associate two lower
semimatrices $\bar{A}=(\bar{a}_{nv})$ and $\hat{A}=(\hat{a}_{nv})$
as follows:
\begin{equation}
\bar{a}_{nv}=\sum_{i=v}^{n}a_{ni},\quad n,v=0,1,...\quad \bar{\Delta}a_{nv}=a_{nv}-a_{n-1,v}, \quad a_{-1,0}=0
\end{equation}
and

\begin{equation}
\hat{a}_{00}=\bar{a}_{00}=a_{00},\quad
\hat{a}_{nv}=\bar{\Delta}\bar{a}_{nv},\quad n=1,2,...
\end{equation}
It may be noted that $\bar{A}$ and $\hat{A}$ are the well-known
matrices of series-to-sequence and series-to-series
transformations, respectively. Then, we have
\begin{equation}\label{eq:13}
A_{n}(s)=\sum_{v=0}^{n}a_{nv}s_{v}=
\sum_{v=0}^{n}\bar{a}_{nv}a_{v}
\end{equation}
and
\begin{equation}
\bar{\Delta}A_{n}(s)=\sum_{v=0}^{n}\hat{a}_{nv}a_{v}.
\end{equation}

\begin{thm} \label{Theomain}
	Let $(p_{n})$ be a sequence of positive numbers such that
	$ P_{n}=O(np_{n}) \ ~~ as \ n\rightarrow\infty$, if $\phi_1(t)\in \mathcal{BV}(0,\pi)$, and the
	sequences $(A_n)$, $(\lambda_n)$, and $(X_n)$ satisfy the conditions
	of Theorem \ref{TheoremBor2}.
	Let $k\geq 1$  and  if $A=(a_{nv})$ is a positive normal matrix such that
	\begin{equation}
	\overline{a}_{n0}=1,\     n=0,1,...,\\
	\end{equation}
	\begin{equation}
	a_{n-1,v}\geq a_{nv},\ \textnormal{for}~~   n\geq v+1,\\
	\end{equation}
	\begin{equation}
	a_{nn}=O(\frac{p_{n}}{P_{n}}),\\
	\end{equation}
	\begin{equation}
	\hat{a}_{n,v+1}=O(v|\overline{\Delta}{a_{nv}}|)\\
	\end{equation}
	then the series $ \sum C_{n}(x)\lambda_{n}$ is summable $\left|A,p_{n}\right|_{k}$, $k\geq 1$.\\
\end{thm}

We need the following lemma for the proof of Theorem \ref{Theomain}.\\
{ Lemma 3} \cite{Bor2} Under the conditions of Theorem 1.1, we have that
\begin{equation}
|\lambda_n|X_n =O(1)\ as \ n\rightarrow\infty,\\
\end{equation}
\begin{equation}
n X_n |A_n| =O(1)\ as \ n\rightarrow\infty,\\
\end{equation}
\begin{equation}
\sum_{n=1}^{\infty}n X_n|\Delta A_n| <\infty.
\end{equation}

\section{Conclusion}
In this paper, the concept of absolute matrix  summability is investigated. In this investigation,  we proved interesting theorem related to $\left| A,p_{n}\right| _{k}$. We also obtain applications to Fourier series.  One may expect this investigation to be a useful tool in the field of analysis in modeling various problems occurring in many areas of science, dynamical systems, computer science, information theory, economical science and biological science.\\
We can apply Theorem \ref{Theomain} to weighted mean $A=(a_{nv})$ is defined as $a_{nv}=\frac{p_{v}}{P_{n}}$ when $0\leq v\leq n$, where $P_{n}=p_{0}+p_{1}+...+p_{n}.$ We have that,
\begin{equation}
\bar{a}_{nv}=\frac{P_{n}-P_{v-1}}{P_{n}} \quad and \quad \hat{a}_{n,v+1}=\frac{p_{n}P_{v}}{P_{n}P_{n-1}}
\end{equation}
(see \cite{mu}-\cite{Sa2} and \cite{4}- \cite{9} for the related bibliography).


\section{ACKNOWLEDGMENTS}
This work was supported by Ahi Evran University Scientific Research Projects Coordination Unit. Project Number: FEF.A3.17.003

\end{document}